\documentstyle[11pt]{article}
\title{On the Deligne-Simpson problem and its weak version}
\author{Vladimir Petrov Kostov\\ \\ \hspace{7cm}
{\sl To Prof. Yu. S. Il'yashenko}} 
\date{}
\newtheorem{tm}{Theorem}
\newtheorem{lm}[tm]{Lemma}
\newtheorem{cor}[tm]{Corollary}
\newtheorem{prop}[tm]{Proposition}
\newtheorem{rem}[tm]{Remark}
\newtheorem{rems}[tm]{Remarks}
\newtheorem{defi}[tm]{Definition}
\newtheorem{ex}[tm]{Example}

\newtheorem{nota}[tm]{Notation}
\newtheorem{conv}[tm]{Convention}
\newtheorem{stat}[tm]{Statement}
\begin{document}
\maketitle 

\begin{abstract}
We consider the {\em Deligne-Simpson problem (DSP) (resp. the weak DSP): 
Give necessary and sufficient conditions upon the choice of the $p+1$ 
conjugacy 
classes $c_j\subset gl(n,{\bf C})$ or $C_j\subset GL(n,{\bf C})$ so that 
there exist irreducible $(p+1)$-tuples (resp. $(p+1)$-tuples with 
trivial centralizers) of matrices $A_j\in c_j$ with zero sum or of 
matrices $M_j\in C_j$ whose product is $I$.} The matrices $A_j$ (resp. $M_j$) 
are interpreted as matrices-residua of Fuchsian linear systems 
(resp. as monodromy matrices of regular linear systems) of 
differential equations with complex time. In the paper we give sufficient 
conditions for solvability of the DSP in the case when one of the matrices is 
with distinct eigenvalues.
\end{abstract}

\section{Introduction}
\subsection{Basic notions and purpose of this paper}

In the present paper we consider the {\em Deligne-Simpson problem (DSP): 
Give necessary and sufficient conditions upon the choice of the $p+1$ 
conjugacy 
classes $c_j\subset gl(n,{\bf C})$ or $C_j\subset GL(n,{\bf C})$ so that 
there exist irreducible $(p+1)$-tuples of matrices $A_j\in c_j$ satisfying 
the condition 
\begin{equation}\label{A_j} 
A_1+\ldots +A_{p+1}=0
\end{equation}
or of matrices} $M_j\in C_j$ satisfying the condition 

\begin{equation}\label{M_j}
M_1\ldots M_{p+1}=I
\end{equation}

\begin{conv}\label{conv1}
In what follows we write ``tuple'' instead of ``$(p+1)$-tuple'' and the 
matrices $A_j$ (resp. $M_j$) are always supposed to satisfy condition 
(\ref{A_j}) (resp. (\ref{M_j})).
\end{conv}
 
The matrices $A_j$ (resp. $M_j$) are interpreted as matrices-residua of 
a Fuchsian system of linear differential equations  
(resp. as monodromy matrices of a regular linear system) on Riemann's sphere; 
see a more detailed description in \cite{Ko1} or 
\cite{Ko2}. 

\begin{rem}
The version with matrices $A_j$ (resp. $M_j$) is called the 
{\em additive} (resp. the {\em multiplicative}) version of the DSP. The 
multiplicative version of the problem was formulated by P.Deligne and 
C.Simpson was the first to obtain results towards its resolution, see 
\cite{Si1} and \cite{Si2}. The additive version is due to the author.
\end{rem}

We presume the necessary condition $\prod \det (C_j)=1$ (resp. 
$\sum$Tr$(c_j)=0$) to hold. In terms of the eigenvalues $\sigma _{k,j}$  
(resp. $\lambda _{k,j}$) of the matrices from $C_j$ (resp. $c_j$) repeated 
with their multiplicities, this condition reads     
$\prod _{k=1}^n\prod _{j=1}^{p+1}\sigma _{k,j}=1$  
(resp. $\sum _{k=1}^n\sum _{j=1}^{p+1}\lambda _{k,j}=0$). 

\begin{defi}\label{genericevs}
An equality $\prod _{j=1}^{p+1}\prod _{k\in \Phi _j}\sigma _{k,j}=1$, resp. 
$\sum _{j=1}^{p+1}\sum _{k\in \Phi _j}\lambda _{k,j}=0$, is called a 
{\em non-genericity relation};  
the sets $\Phi _j$ contain one and the same number $N<n$ of indices  
for all $j$ (when wishing to specify $N$ we say ``$N$-relation'' instead of 
``non-genericity relation''). Eigenvalues satisfying none of these 
relations are called 
{\em generic}. 
\end{defi}

\begin{rems}\label{reducible}
1) Reducible tuples of matrices $A_j$ or $M_j$ 
exist only for non-generic eigenvalues (the 
eigenvalues of each diagonal block of a block upper-triangular  
tuple satisfy some non-genericity relation). Therefore for generic 
eigenvalues existence of tuples implies automatically their irreducibility.
This is not true for non-generic eigenvalues.

2) It is clear that the presence of a non-genericity relation with $N=N_0$ 
implies the presence of one with $N=n-N_0$ (just replace the sets $\Phi _j$ 
by their complements in $\{ 1,2,\ldots ,n\}$). Therefore in what follows we 
consider only non-genericity relations with $N\leq n/2$. 
\end{rems}
 
Part 1) of the above remarks explains why for non-generic 
eigenvalues it is reasonable 
to require instead of irreducibility of the tuple only triviality of its 
centralizer (i.e. only scalar matrices to commute with all matrices from 
the tuple). This is the {\em weak version of the DSP} (or just 
the {\em weak DSP} for short). 

\begin{defi}
We say that the DSP (resp. the weak DSP) is {\em solvable} for a given 
tuple of conjugacy classes $c_j$ or $C_j$ if there exist irreducible 
tuples of matrices $A_j\in c_j$ 
or $M_j\in C_j$ (resp. if there exist tuples of such matrices with trivial 
centralizers).
\end{defi} 

We assume throughout the paper that there holds 

\begin{conv}\label{conv2}
The conjugacy classes $c_1$ and $C_1$ are with distinct eigenvalues.
\end{conv}

The purpose of the present paper is to show as precisely as possible 
where passes the border between the cases 
when the DSP is solvable and when it is not but the weak DSP is solvable.

\subsection{The known results}

\begin{defi}\label{JNform}
Call {\em Jordan normal form (JNF) of size $n$} a family 
$J^n=\{ b_{i,l}\}$ ($i\in I_l$, $I_l=\{ 1,\ldots ,s_l\}$, $l\in L$) of 
positive integers $b_{i,l}$ 
whose sum is $n$. Here $L$ is the set of indices of 
eigenvalues (all distinct) and 
$I_l$ is the set of indices of 
Jordan blocks with eigenvalue $l$, $b_{i,l}$ is the 
size of the $i$-th block with this eigenvalue. An $n\times n$-matrix 
$Y$ has the JNF $J^n$ (notation: $J(Y)=J^n$) if to its distinct 
eigenvalues $\lambda _l$, $l\in L$, there belong Jordan blocks of sizes 
$b_{i,l}$. We use the following notation (illustrated by an example): 
the JNF $\{ \{ 3,2\} ,\{ 7,6,1\} \}$ is the one with two eigenvalues to the 
first (to the second) of which there belong two blocks, of sizes $3$ and $2$ 
(resp. three blocks, of sizes $7$, $6$ and $1$). 
\end{defi}

\begin{nota}\label{dr}
1) We denote by $C(Y)$ the conjugacy class (in $gl(n,{\bf C})$ or 
$GL(n,{\bf C})$) of the matrix $Y$. We set $C(Y)=C(X)\times C(Z)$ if 
$Y=\left( \begin{array}{cc}X&0\\0&Z\end{array}\right)$ (here 
$X$ is $l\times l$ and $Z$ is $(n-l)\times (n-l)$. 

2) For a conjugacy class $C$ in $GL(n,{\bf C})$ or $gl(n,{\bf C})$ denote by 
$d(C)$ its dimension and by $J(C)$ the JNF it defines. 
For a matrix $Y\in C$ set 
$r(C):=\min _{\lambda \in {\bf C}}{\rm rank}(Y-\lambda I)$. The integer 
$n-r(C)$ is the maximal number of Jordan blocks of $J(Y)$ with one and the 
same eigenvalue. Set $d_j:=d(C_j)$ (resp. $d(c_j)$), $r_j:=r(C_j)$ 
(resp. $r(c_j)$). The quantities 
$r(C)$ and $d(C)$ depend only on the JNF $J(Y)=J^n$, not 
on the eigenvalues, so we write sometimes $r(J^n)$ and $d(J^n)$. 
\end{nota}

\begin{prop}\label{d_jr_j}
(C. Simpson, see \cite{Si1}.) The 
following couple of inequalities is a necessary condition for the existence 
of irreducible $(p+1)$-tuples satisfying (\ref{M_j}) or (\ref{A_j}):

\[ \begin{array}{lll}d_1+\ldots +d_{p+1}\geq 2n^2-2&~~~~~~~~&(\alpha _n)\\ \\ 
{\rm for~all~}j,~r_1+\ldots +\hat{r}_j+\ldots +r_{p+1}\geq n&~~~~~~~~&
(\beta _n)\end{array}\]
\end{prop}

The above proposition holds without Convention~\ref{conv2}. When 
Convention~\ref{conv2} holds, then $r_1=n-1$ and 
condition $(\beta _n)$ is tantamount to 
$r_2+\ldots +r_{p+1}\geq n$.

\begin{defi}\label{indexrig}
The quantity $\kappa =2n^2-d_1-\ldots -d_{p+1}$ (see Notation~\ref{dr}) 
is called the {\em index of rigidity} of a given tuple of conjugacy classes 
or of JNFs. It has been introduced by N.Katz, see \cite{Ka}. If condition 
$(\alpha _n)$ holds, then $\kappa$ can take the values 
$2$, $0$, $-2$, $-4$, $\ldots$. The case $\kappa =2$ is 
called the {\em rigid} one.
\end{defi}

\begin{defi}
A {\em multiplicity vector (MV)} is a vector whose components are non-negative 
integers whose sum is $n$. Further in the text components of the MVs are 
the multiplicities of the eigenvalues of $n\times n$-matrices. 
%For such a MV 
%$(m_1,\ldots ,m_s)$, $m_1+\ldots +m_s=n$, we assume that 
%$m_1\geq \ldots \geq m_s$.
\end{defi}

\begin{rem}\label{d_jdiag}
For a diagonalizable conjugacy class $C$ 
with MV equal to $(m_1,\ldots ,m_s)$ 
one has $d(C)=n^2-m_1^2-\ldots -m_s^2$.
\end{rem}

\begin{defi}\label{corrdefi}
For a given JNF $J^n=\{ b_{i,l}\}$ define its {\em corresponding} 
diagonal JNF ${J'}^n$. A diagonal JNF is  
a partition of $n$ defined by the multiplicities of the eigenvalues. 
For each $l$ $\{ b_{i,l}\}$ is a partition of $\sum _{i\in I_l}b_{i,l}$ and 
$J'$ is the disjoint sum of the dual partitions. 
Thus if for each fixed $l$ one has 
$b_{1,l}\geq$$\ldots$$\geq b_{s_l,l}$,  
then the eigenvalue 
$l\in L$ is replaced by $b_{1,l}$ new eigenvalues $h_{1,l}$, 
$\ldots$, $h_{b_{1,l},l}$ (hence, ${J'}^n$ has 
$\sum _{l\in L}b_{1,l}$ distinct eigenvalues). 
\end{defi}

\begin{rems}\label{corrrems}
One has the following properties of corresponding JNFs (see \cite{Ko2}:)

1) For $l$ fixed, set $g_k$ 
for the multiplicity of the eigenvalue $h_{k,l}$. Then the  
first $b_{s_l,l}$ numbers $g_k$ equal 
$s_l$, the next $b_{s_{l-1},l}-b_{s_l,l}$ equal $s_l-1$, $\ldots$, 
the last $b_{1,l}-b_{2,l}$ equal 1. 

2) There hold the equalities $r(J^n)=r({J'}^n)$ and $d(J^n)=d({J'}^n)$.

3) To each diagonal JNF there corresponds a unique JNF with a single 
eigenvalue.
\end{rems}

\begin{lm}\label{betaequalalpha}
Given the $p+1$ diagonalizable conjugacy classes $c_j$ or $C_j$ satisfying 
condition $(\beta _n)$ and Convention~\ref{conv2}, 
condition $(\alpha _n)$ does not hold for them only in  

{\bf Case A) :}  $p=2$, $n\geq 4$ is even and the MVs of 
$c_2$ and $c_3$ 
(resp. of $C_2$ and $C_3$) both equal $(n/2,n/2)$.
\end{lm}

The lemma is proved at the end of the subsection.

\begin{rem}\label{betaequalalpharem}
Making use of Definition~\ref{corrdefi} and Remarks~\ref{corrrems} one can 
extend the lemma to the case of not necessarily diagonalizable matrices 
(except $A_1$ or $M_1$). In such a context, in Case A) each conjugacy class 
$c_2$, $c_3$ or $C_2$, $C_3$ is either diagonalizable and as in the lemma 
or with a single eigenvalue and $n/2$ Jordan blocks of size $2$ 
belonging to it. Indeed, this is the only non-diagonal JNF corresponding to 
the one with two eigenvalues each of multiplicity $n/2$.
\end{rem}  

The first important result in the resolution of the DSP was the following 

\begin{tm}\label{Simpson}
(C.Simpson, see \cite{Si1}) For generic eigenvalues and under 
Convention~\ref{conv2} conditions $(\alpha _n)$ and $(\beta _n)$ together are 
necessary and sufficient for the solvability of the DSP for given conjugacy 
classes $C_j$.
\end{tm}

The same result for classes $c_j$ is proved in \cite{Ko4}, Theorem 19. 
For arbitrary eigenvalues there holds the following theorem 
(see \cite{Ko3}, Theorem 6). 

\begin{tm}\label{weakDSPdiag}
Under 
Convention~\ref{conv2} conditions $(\alpha _n)$ and $(\beta _n)$ together are 
necessary and sufficient for the solvability of the weak DSP for given 
conjugacy classes $c_j$ or $C_j$.
\end{tm}

\begin{rems}\label{4rigidcases}
1) In \cite{Si1} C.Simpson has considered the rigid case for diagonalizable 
matrices and under 
Convention~\ref{conv2}. He has shown that conditions $(\alpha _n)$ and 
$(\beta _n)$ together hold only if $p=2$ and the MVs of the three matrices 
correspond to one of the four cases:

\[ \begin{array}{lllll}(1,\ldots ,1)&(1,\ldots ,1)&(n-1,1)&&
{\rm hypergeometric~family}\\ \\ 
(1,\ldots ,1)&(\frac{n}{2},\frac{n}{2}-1,1)&(\frac{n}{2},\frac{n}{2})&&
{\rm even~family}\\ \\ 
(1,\ldots ,1)&(\frac{n-1}{2},\frac{n-1}{2},1)&(\frac{n+1}{2},\frac{n-1}{2})&&
{\rm odd~family}\\ \\ 
(1,1,1,1,1,1,)&(2,2,2)&(4,2)&&{\rm extra~case}\end{array}\]

Observe that in all four cases one has $r_2+r_3=n$, i.e. there is an equality 
in condition $(\beta _n)$. Although C.Simpson considers only matrices $M_j$, 
the result is automatically extended to the case of matrices $A_j$.

2) If one wants to get rid of the condition the matrices to be diagonalizable 
(except $A_1$ or $M_1$), then to the above list one should add all cases 
when a diagonal JNF from the list 
is replaced by a JNF corresponding to it. All JNFs 
corresponding to the one with $n$ distinct eigenvalues are the ones in 
which to each eigenvalue there belongs a single Jordan block. Using the 
notation from Definition~\ref{JNform}, give the list of all 
JNFs corresponding to 
the other diagonal JNFs (defined by the MVs) 
encountered in part 1) of the present remarks:

\[ \begin{array}{lllll}(n-1,1)&&\{ 2,1,\ldots ,1\} &&\\ \\ 
(\frac{n}{2},\frac{n}{2}-1,1)&&\{ 3,2,\ldots ,2,1\} &{\rm or}&
\{ \{ 1,\ldots ,1\}\{ 2,1,\ldots ,1\} \} ~
{\rm (}n/2~{\rm and~}n/2-2~{\rm units)}\\ \\ 
(\frac{n}{2},\frac{n}{2})&&\{ 2,\ldots 2\} &&\\ \\ 
(\frac{n-1}{2},\frac{n-1}{2},1)&&\{ 3,2,\ldots ,2\} &{\rm or}& 
\{ \{ 1,\ldots ,1\}\{ 2,1,\ldots ,1\} \} ~{\rm (}(n-1)/2~
{\rm and}~(n-3)/2~{\rm units)}\\ \\ 
(\frac{n+1}{2},\frac{n-1}{2})&&\{ 2,\ldots ,2,1\} &&\\ \\ 
(2,2,2)&& \{ 3,3\} &{\rm or}&\{ \{ 2,2\} \{ 1,1\} \} \\ \\ 
(4,2)&&\{ 2,2,1,1\} &&\end{array}\] 
\end{rems}

{\em Proof of Lemma~\ref{betaequalalpha}:}

$1^0$. Suppose first that one has 

\[ r_j\leq n/2~~{\rm for~~}j=2,\ldots ,p+1~~~~~~~~~~~(*)\] 
Then one has 
$d_j\geq 2r_j(n-r_j)$ and there is equality if and only if the MV of $c_j$ 
or $C_j$ equals $(r_j,n-r_j)$. This follows from Remark~\ref{d_jdiag}. 

For $r_2+\ldots +r_{p+1}$ fixed the sum $d_2+\ldots +d_{p+1}$ is 
minimal for 
$r_2=r_3=[n/2]$ where $[.]$ stands for the entire part of. Indeed, one 
has $d_2+\ldots +d_{p+1}=(r_2+\ldots +r_{p+1})n-r_2^2-\ldots -r_{p+1}^2$ and 
one has to maximize $r_2^2+\ldots +r_{p+1}^2$ for $r_2+\ldots +r_{p+1}$ fixed 
while respecting condition $(*)$.

If $n$ is even and $r_2=r_3=n/2$, $r_j=0$ for $j>3$, 
then condition $(\alpha _n)$ fails if and only if $n\geq 4$ (this is Case A)); 
if $r_4\neq 0$, then condition 
$(\alpha _n)$ holds. If $n$ is odd, then the sum $d_2+\ldots +d_{p+1}$ is 
minimal for $r_2=r_3=[n/2]$, $r_4=1$ and condition $(\alpha _n)$ holds. One 
cannot have $r_j=0$ for all $j>3$ because then condition $(\beta _n)$ does not 
hold.

$2^0$. Suppose that $r_2>n/2$. Denote the MV of the class $c_2$ or $C_2$ 
by $(m_1,\ldots ,m_s)$, with $m_1\geq \ldots \geq m_s$. 
Then $d_2$ is minimal if 
$m_1=m_2=\ldots =m_{s-1}=n-r_2$, see Remark~\ref{d_jdiag}. The sum 
$d_3+\ldots +d_{p+1}$ is minimal if $r_3=m_1=n-r_2$, $r_4=\ldots =r_{p+1}=0$ 
and the MV defining the class $c_3$ or $C_3$ equals $(r_2,n-r_2)$. 

Set $n=(s-1)m_1+m_s$. Recall that $1\leq m_s\leq m_1$. Hence, 

\[ d_1=n^2-n~,~d_2=n^2-(s-1)m_1^2-m_s^2\geq n^2-m_1n~,~d_3=2m_1(n-m_1)~~
{\rm and}\]  

\[ d_1+d_2+d_3\geq 2n^2-n+m_1n-2m_1^2\geq 2n^2-n+n-2=
2n^2-2~~{\rm because}~~1\leq m_1<n/2 \]

The lemma is proved.~~~~~$\Box$

\subsection{The new results}

\begin{defi}
The eigenvalues of the matrices $A_j$ or $M_j$ are called $k$-{\em generic}, 
$k\in {\bf N}$, if they satisfy non-genericity relations only with 
$N\geq k$, see Definition~\ref{genericevs} and part 2) of 
Remarks~\ref{reducible}.
\end{defi} 

\begin{tm}\label{123generic}
Under Convention~\ref{conv2}, if the eigenvalues are $2$-generic, and 
if $\kappa \leq 0$ (see Definition~\ref{indexrig}), then 
conditions $(\alpha _n)$ and $(\beta _n)$ are necessary and sufficient for 
the solvability of the DSP. 
\end{tm}

The theorem is proved in Section~\ref{proofof123generic}. 
Examples~\ref{exTn+1} and \ref{ex23} below show that the theorem 
cannot be made stronger.

\begin{tm}\label{n+1}
Under Convention~\ref{conv2} and for arbitrary eigenvalues, 
if $r_2+\ldots +r_{p+1}\geq n+1$, then the DSP is solvable for such 
conjugacy classes.
\end{tm} 

The theorem is proved in Section~\ref{proofofn+1}.  
Example~\ref{exTn+1} below  
shows that for $r_2+\ldots +r_{p+1}=n$ Theorem~\ref{n+1} is no longer true.

\begin{rem}\label{trivc}
The above two theorems imply that under Convention~\ref{conv2} 
the weak DSP is 
solvable but the DSP is not only if $r_2+\ldots +r_{p+1}=n$ and  
either $\kappa =2$ or the eigenvalues satisfy a $1$-relation. 
\end{rem}

\begin{cor}\label{cor1}
Under Convention~\ref{conv2} a block upper-triangular tuple of 
diagonalizable matrices $A_j$ or $M_j$ 
with $3$-generic eigenvalues can be deformed into one from the same 
conjugacy classes and with trivial centralizer.
\end{cor}

Indeed, $3$-genericity 
implies that for each diagonal block (say, of size $s\geq 3$) 
there holds condition 
$(\beta _s)$ and Case A) from Lemma~\ref{betaequalalpha} 
is avoided; hence, condition $(\beta _n)$ holds for the tuple of 
conjugacy classes (the quantity $r$ computed for the whole matrix 
is not smaller than the sum of the quantities $r$ computed for the diagonal 
blocks), and Case A)  
is avoided  
(because the blocks are of size $\geq 3$ -- we leave the details 
for the reader). Hence, for the given tuple of conjugacy classes there hold 
conditions $(\alpha _n)$ and $(\beta _n)$ (see Lemma~\ref{betaequalalpha}). 
The claim follows now from Lemma~24 from \cite{Ko3}.~~~~~$\Box$ 

\begin{cor}
Under Convention~\ref{conv2}, if the eigenvalues are $2$-generic, 
and if Case A) is avoided, then for such a 
block upper-triangular tuple of diagonalizable matrices $A_j$ or $M_j$ 
there hold conditions $(\beta _n)$ and $(\alpha _n)$. Moreover, the tuple 
can be deformed into one from the same 
conjugacy classes and with trivial centralizer.
\end{cor}

The first claim is proved as Corollary~\ref{cor1}, the second 
follows from Lemma~24 from \cite{Ko3}.

\begin{nota}\label{semidirect}
For a tuple of  
matrices $A_j$ or $M_j$ in block upper-triangular form 
$\left( \begin{array}{cc}P_j&Q_j\\0&R_j\end{array}\right)$ (where 
$P_j\in gl(l,{\bf C})$, $R_j\in gl(n-l,{\bf C})$) set 
$d_j^1=d(P_j)$, $r^1_j=r(P_j)$, $d_j^2=d(R_j)$, $r^2_j=r(R_j)$, 
$s_j=$dim${\cal X}_j$ where 
${\cal X}_j=\{ Z\in M_{l,n-l}|Z=P_jX_j-X_jR_j, X_j\in M_{l,n-l}\}$. 
Denote by ${\cal P}$, ${\cal R}$ the representations defined by the tuples 
of matrices $P_j$, $R_j$.
\end{nota}

\begin{rem}\label{semidirectrem}
If the MVs of the diagonalizable matrices $P_j$ and $R_j$ equal respectively 
$(m_1',\ldots ,m_s')$, $(m_1'',\ldots ,m_s'')$ (there might be zeros among 
these numbers as some eigenvalue might be absent in $P_j$ or $R_j$), 
then $s_j=l(n-l)-\sum _{i=1}^sm_i'm_i''$. This implies that if one exchanges 
the positions of the blocks $P_j$ and $R_j$, then the quantities $s_j$ do not 
change. 
\end{rem}

\begin{lm}\label{Ext}
If the representations ${\cal P}$ and ${\cal R}$ are with trivial 
centralizers, then one has 

\[ \delta :={\rm dim\,Ext}^1({\cal P},{\cal R})=s_1+\ldots +s_{p+1}-2l(n-l)~.\]
\end{lm}

{\em Proof:}

Notice first that ${\cal X}_j$ is the space of right upper blocks of matrices 
of the form 

\[ \left( \begin{array}{cc}I&X_j\\0&I\end{array}\right) ^{-1}
\left( \begin{array}{cc}P_j&0\\0&R_j\end{array}\right) 
\left( \begin{array}{cc}I&X_j\\0&I\end{array}\right) ~.\]
To obtain $\delta$ one must first subtract $l(n-l)$ 
from $\sum _{j=1}^{p+1}$dim${\cal X}_j$ (because the sum of these right upper 
blocks must be $0$) and then again subtract $l(n-l)$ (to factor out the 
simultaneous conjugation with matrices 
$\left( \begin{array}{cc}I&X\\0&I\end{array}\right)$; as $A_1$ or $M_1$ is 
with distinct eigenvalues, no such matrix with $X\neq 0$ commutes with all 
matrices from the tuple).~~~~~$\Box$

\begin{ex}\label{exTn+1}
Consider under Convention~\ref{conv2} a tuple of diagonalizable 
conjugacy classes $c_j$ 
for which $r_2+\ldots +r_{p+1}=n$, $n>2$. Denote by $\mu _1$ an 
eigenvalue of $c_1$ and by $\mu _2$, $\ldots$, $\mu _{p+1}$ eigenvalues of 
$c_2$, $\ldots$, $c_{p+1}$ of maximal possible multiplicity; we assume these 
multiplicities to be $>n/2$. Suppose that 
the eigenvalues of the classes $c_j$ satisfy the only non-genericity 
relation $\mu _1 +\ldots +\mu _{p+1}=0$. 

Denote by $c_j'\subset gl(n-1,{\bf C})$ the conjugacy classes obtained from 
$c_j$ by deleting the eigenvalues $\mu _j$. Hence, condition 
$(\beta _{n-1})$ holds for the classes $c_j'$ and the sum of their eigenvalues 
is $0$. Moreover, the classes $c_j'$ do not correspond to Case A) from 
Lemma~\ref{betaequalalpha} (we let the reader check this oneself).

Hence, there exist block upper-triangular matrices 
$A_j=\left( \begin{array}{cc}A_j'&D_j\\0&\mu _j\end{array}\right)$, 
$A_j'\in c_j'$, whose 
tuple defines a semi-direct sum (but not a direct one); the matrices $A_j'$ 
define an irreducible representation. Indeed, one checks directly that 
dim\,Ext$^1(A,M)=1$ (this results from $r_2+\ldots +r_{p+1}=n$). The same 
equality shows that the variety ${\cal V}$ 
consisting of tuples of matrices $A_j\in c_j$ which are block 
upper-triangular up to conjugacy (i.e. like $A_j$ above) is of dimension 
dim${\cal W}$ where ${\cal W}$ is 
the variety of tuples with trivial centralizers from the classes 
$c_j$. 

This means that there exist no irreducible tuples from the classes $c_j$. 
Indeed, should they exist, their variety (which is part of ${\cal W}$) 
should contain in its closure the 
variaty ${\cal V}$ (see Theorem~6 from \cite{Ko3}), 
hence, one would have dim${\cal V}<$dim${\cal W}$ which 
is a contradiction.

The example shows that Theorem~\ref{123generic} is not 
true without the condition the eigenvalues to be $2$-generic and that 
Theorem~\ref{n+1} is not true if there is an equality in $(\beta _n)$. 

A similar example can be given for matrices $M_j$.
\end{ex}

\begin{ex}\label{ex23}
There exist triples of diagonalizable $2\times 2$-matrices $M_j^1$ 
(resp. $M_j^2$) with (generic) eigenvalues equal to 
$(a,b)$, $(\mu ,\nu )$, $(\eta ,\xi )$ (resp. to $(c,d)$, $(\mu ,\nu )$, 
$(\eta ,\zeta )$); same (different) letters denote same (different) 
eigenvalues. 

Then there exists a block upper-triangular triple of matrices 
$M_j=\left( \begin{array}{cc}M_j^1&B_j\\0&M_j^2\end{array}\right)$ defining a 
semi-direct sum of the representations ${\cal P}^1$ and ${\cal P}^2$ 
defined by the matrices $M_j^1$ and 
$M_j^2$ (because dim\, Ext$^1({\cal P}^1,{\cal P}^2)=1$). 

One checks directly that 

a) the centralizer of the matrices $M_j$ is trivial; 

b) their eigenvalues can be chosen $2$-generic (we assume that they satisfy 
only the following non-genericity relations: $ab\mu \nu \eta \xi =1$ and 
$cd\mu \nu \eta \zeta =1$);

c) one has $\kappa =2$ 
for the triple of conjugacy classes of the matrices $M_j$. 

As $\kappa =2$, one cannot have coexistence of irreducible and reducible 
triples, see \cite{Ka}. This means that the DSP is not solvable for the 
triple of conjugacy classes of the matrices $M_j$ (but the weak DSP is, see 
a)). Hence, Theorem~\ref{123generic} is not true for $\kappa =2$.

A similar example can be given for matrices $A_j$. 
\end{ex}

\section{Proof of Theorem~\protect\ref{123generic}
\protect\label{proofof123generic}}

\subsection{The method of proof} 

$1^0$. Suppose that for the conjugacy classes 
$c_j$ or $C_j$ (with $2$-generic eigenvalues) 
there hold conditions $(\alpha _n)$ and $(\beta _n)$. 
The variety of matrices $A_j\in c_j$ (satisfying (\ref{A_j})) or 
of matrices $M_j\in C_j$ (satisfying (\ref{M_j})) is of dimension 
$d':=d_1+\ldots +d_{p+1}-n^2+1$ at each tuple with trivial centralizer, 
see \cite{Ko6}, Proposition~2. 

Given a reducible tuple of matrices from these conjugacy classes 
(block upper-triangular up to conjugacy, with trivial 
centralizer, with given 
sizes of the diagonal blocks and with given conjugacy classes of the 
restrictions of the matrices to the diagonal blocks) we compute the 
dimension $d''$ of the variety of 
such tuples and we show that $d''<d'$. If this is the case of 
all such reducible tuples, then the variety of tuples with 
trivial centralizers 
must contain irreducible tuples as well. Hence, the DSP is solvable for the 
given conjugacy classes.

\begin{lm}\label{lmdim}
Under Convention~\ref{conv2}, suppose that the tuple of diagonalizable 
matrices $A_j$ 
or $M_j$ is as in Notation~\ref{semidirect}, and that 
the representations ${\cal P}$ 
and ${\cal R}$  
are with trivial centralizers.
%, and suppose that
%
%\[ J\left( \left( \begin{array}{cc}P_j&Q_j\\0&R_j\end{array}\right) \right) =
%J\left( \left( \begin{array}{cc}P_j&0\\0&R_j\end{array}\right) \right) ~.\] 
If $\delta :=$dim\,Ext$^1({\cal P}^1,{\cal P}^2)>1$, then 
$d''<d'$.
\end{lm}

All lemmas from the proof of the theorem are proved in Subsection~\ref{prlm}.

\begin{cor}\label{existirred}
If the representations ${\cal P}$ and ${\cal R}$ from the lemma are 
irreducible, then there exist irreducible tuples from the conjugacy 
classes $c(P_j)\times c(R_j)$.
\end{cor}

The corollary is immediate.

We prove the theorem for diagonalizable matrices in $2^0$ -- $5^0$ 
and then we treat the general case in $6^0$ -- $11^0$.

\subsection{The proof for diagonalizable matrices}

$2^0$. Prove the theorem for diagonalizable matrices. 

\begin{lm}\label{deform}
Suppose that the tuples of diagonalizable matrices 
$P_j\in gl(l,{\bf C})$ and $R_j\in gl(n-l,{\bf C})$ 
(resp. $P_j\in GL(l,{\bf C})$ and $R_j\in GL(n-l,{\bf C})$) 
are with trivial centralizers, $P_1$ and $R_1$ being each with 
distinct eigenvalues and with 
no eigenvalue in common, 
and that $l\geq n-l\geq 2$. Then $\delta \geq 2$ with the exception of the 
cases listed below\footnote{When listing the cases we begin with B, 
not with A, in order to avoid mixing up with Case A) from 
Lemma~\protect\ref{betaequalalpha}}. In all of them one has $p=2$. 
(We give the list of the eigenvalues of the 
matrices $P_2$, $R_2$ and $P_3$, $R_3$, equal (different) letters 
denote equal (different) eigenvalues if they 
correspond to one and the same index $j$. In Cases C) -- F) one can 
exchange the roles of $P_2$, $R_2$ and $P_3$, $R_3$.)

\[ \begin{array}{lllll}{\rm Case~B)}&l=n-l=2&
(a,b)&~~&(c,d)\\&&(a,b)&&(c,d)\\ \\ 
{\rm Case~C)}&l=n-l=2&(a,b)&~~&(c,d)\\&&(a,g)&&(c,d)\\ \\ 
%{\rm Case~D)}&l=n-l=2&(a,b)&~~&(c,d)\\&&(a,b)&&(c,h)\\ \\ 
%{\rm Case~E)}&l=3~,~n-l=2&(a,b,c)&(a,b)&~~&(f,g,g)&(f,g)\\
{\rm Case~D)}&l=n-l=3&(a,b,c)&~~&(f,g,g)\\&&(a,b,c)&&(f,g,g)\\ \\ 
{\rm Case~E)}&l=2q+1,~n-l=2&(\underbrace{a,\ldots ,a}_{q~{\rm times}},
\underbrace{b,\ldots ,b}_{q~{\rm times}},c)&~~&
(\underbrace{f,\ldots ,f}_{q+1~{\rm times}},
\underbrace{g,\ldots ,g}_{q~{\rm times}})\\&&(a,b)&&(f,g)\\ \\ 
{\rm Case~F)}&l=2q,~n-l=2&(\underbrace{a,\ldots ,a}_{q~{\rm times}},
\underbrace{b,\ldots ,b}_{q-1~{\rm times}},c)&~~&
(\underbrace{f,\ldots ,f}_{q~{\rm times}},
\underbrace{g,\ldots ,g}_{q~{\rm times}})\\&&(a,b)&&(f,g)
\end{array}\] 

In Case B) condition $(\alpha _n)$ does not hold for the conjugacy 
classes $C(P_j)\times C(R_j)$, in the other cases it holds and is an 
equality. One has $\delta =0$ in Case B) and $\delta =1$ in Cases C) -- F).
\end{lm}

\begin{cor}
In the conditions of the lemma and if the representations ${\cal P}$ and 
${\cal R}$ are irreducible the DSP is solvable for the tuple of conjugacy 
classes $C(S_j)=C(P_j)\times C(R_j)$ (except for Cases B) -- F)).
\end{cor}

{\em Proof:} 

The condition $\delta >0$ implies that there exists a semi-direct sum of the 
representations ${\cal P}$ and ${\cal R}$ (we use Notation~\ref{semidirect} 
here) 
which is not reduced to a direct one. The centralizer of this semi-direct sum 
is trivial. Indeed, one can assume that $P_1$ and $R_1$ are diagonal, 
so a matrix $X$ from the centralizer must be also diagonal. The $P$-block of 
$X$ commutes with all matrices $P_j$, hence, it is scalar (because the 
centralizer of ${\cal P}$ is trivial). 
In the same way the $R$-block of $X$ must be scalar. Finally, these blocks 
must be equal, otherwise the commutation relations imply that 
all blocks $Q_j$ must be $0$ which contradicts the 
sum of ${\cal P}$ and ${\cal R}$ not to be a direct one.

Hence, the variety ${\cal V}$ 
of tuples of matrices defining semi-direct sums of 
${\cal P}$ and ${\cal R}$ is non-empty and its dimension is smaller than the 
dimension of the variety ${\cal W}\supset {\cal V}$ 
of tuples with trivial centralizers 
of matrices from the classes $C(S_j)$ (see Lemma~\ref{lmdim}). Hence, 
${\cal V}$ is locally a proper subvariety of ${\cal W}$ and 
a tuple from ${\cal V}$ can be deformed into a tuple from 
${\cal W}\backslash {\cal V}$ (see Theorem~6 from \cite{Ko3}). 
The latter must be irreducible. Indeed, 
${\cal V}$ contains locally all reducible tuples because 
${\cal P}$ and ${\cal R}$ are irreducible.~~~~~$\Box$\\

$3^0$. Deduce the theorem from the corollary. 
The weak DSP is solvable 
for conjugacy classes in the conditions of the theorem. Indeed, 
$2$-genericity implies that a tuple from the given conjugacy classes is 
(up to conjugacy) block upper-triangular with diagonal blocks all of sizes 
$\geq 2$ and defining irreducible representations. (We assume that there is 
more than one diagonal block, otherwise the tuple is irreducible and 
there is nothing to prove.) 

The restriction of the tuple to the union of diagonal blocks is a tuple 
from the same conjugacy classes (because the conjugacy classes are 
diagonalizable). Consider a couple of consecutive diagonal blocks. 
(We denote the restrictions of the matrices $A_j$ or $M_j$ to these two blocks 
by $A_j^i$, $M_j^i$, $i=1,2$.) They are 
both of size $\geq 2$, and if one is not in one of the Cases B) -- F), 
then one can apply the above corollary and obtain the 
existence of irreducible tuples of matrices from the conjugacy classes 
$C(A_j^1)\times C(A_j^2)$ (resp. $C(M_j^1)\times C(M_j^2)$). Thus we obtain 
a block-diagonal tuple of $n\times n$-matrices with one diagonal block less. 
Continuing like this we end with an irreducible tuple of matrices which 
solves the DSP for the conjugacy classes $c_j$ or $C_j$.   

$4^0$. There might be a problem, however, with Cases B) -- F). First of all 
notice that this does not happen if $p\geq 3$. Indeed, in this case one can 
always choose two diagonal blocks defining irreducible representations and 
in which at least four conjugacy classes 
$C(A_j^1)\times C(A_j^2)$ (resp. $C(M_j^1)\times C(M_j^2)$) are not scalar 
(including $j=1$). So one can permute the diagonal blocks (to get two 
consecutive blocks not from Cases B) -- F)) and the proof is carried out as 
in $3^0$.

$5^0$. So suppose that $p=2$. We start again with the 
restriction of the tuple to the set of diagonal blocks defining irreducible 
representations. It is not possible to have all couples of diagonal blocks 
to correspond to Case B) from the lemma because this will mean that the 
classes $c_j$ or $C_j$ are from Case A) of Lemma~\ref{betaequalalpha}. 
So choose a couple of consecutive 
diagonal blocks which are 
not from Case B) and replace them by a single block $B$ 
defining a semi-direct sum of the representations which they define while 
keeping the other diagonal blocks the same. This is possible because for 
the chosen blocks one has $\delta \geq 1$, see the lemma. 

At each next step one has a block-diagonal tuple with diagonal blocks defining 
irreducible representations except $B$ which defines one with trivial 
centralizer. At each step choose a block $W$ different from $B$ and 
next to $B$ (hence, 
their couple is not from Case B) because $B$ is of size $>2$), 
so one can replace it by a new block 
(which is the new block $B$) 
defining a semi-direct sum of the representations they define. So at each step 
the blocks $B$, $W$ are not from Case B).

At the last step we obtain a representation with trivial centralizer.  
The last couple of blocks $B$, $W$ is not from Cases B) -- F). Indeed, should 
it be from these cases, then for the conjugacy classes $c_j$ or $C_j$ one 
should have $\kappa \geq 2$ (to be checked directly). 

Hence, for the last couple of blocks $B$, $W$ one has $\delta \geq 2$. 
This means that $d''<d'$, see $1^0$. This proves the theorem in the case 
of diagonalizable matrices.\\

\subsection{The proof in the general case}

$6^0$.

\begin{conv} 
From here till the end of this subsection when Case A) 
of Lemma~\ref{betaequalalpha} or Cases B) -- F) of Lemma~\ref{deform} 
are cited the JNFs of the matrices $A_j$ or $M_j$ ($j\geq 2$) will be assumed 
either to be the ones given in these two lemmas or to correspond to them, 
see Remarks~\ref{betaequalalpharem} and \ref{4rigidcases}.
\end{conv}

Such a change of the definition of these cases does not 
change the quantity $\delta$, 
see part 2) of Remarks~\ref{corrrems}. 
Hence, Lemma~\ref{deform} is applicable after the change as well.  

$7^0$. Consider a tuple in block upper-triangular form whose diagonal 
blocks define irreducible representations. Consider the restriction of the 
tuple to the set of diagonal blocks. The conjugacy class $c_j'$ (resp. $C_j'$) 
of the restriction 
of the matrix $A_j$ (resp. $M_j$) from the tuple to the set of diagonal 
blocks belongs to the closure of $c_j$ (resp. of $C_j$) but is not necessarily 
equal to it (one might obtain a ``less generic'' Jordan structure when 
cutting off the blocks above the diagonal; the eigenvalues and their 
multiplicities do not change). If for the conjugacy classes $c_j'$ or $C_j'$ 
the index of rigidity is $\leq 0$, then as in the case of 
diagonalizable conjugacy classes one shows that the DSP is solvable for the 
classes $c_j'$ or $C_j'$. This implies its solvability for the classes 
$c_j$ (resp. $C_j$) (which can be proved by analogy with part~2 of Lemma~53 
from \cite{Ko2}).

$8^0$. Suppose (in $8^0$ -- $11^0$) 
that the index of rigidity of the tuple of conjugacy 
classes $c_j'$ or $C_j'$ is $>0$. Then for some $j_0>1$ there exists a 
conjugacy class $c_{j_0}''$ (or $C_{j_0}''$; we write further only  
$c_{j_0}''$ for short) such that 

1) $c_{j_0}'$ belongs to the closure of $c_{j_0}''$; 

2) $c_{j_0}''$ is obtained from $c_{j_0}'$ 
when a couple of Jordan blocks with one and the same eigenvalue, of sizes 
$l,s$, $l\geq s$, are replaced by Jordan blocks (with the same eigenvalues) of 
sizes $l+1,s-1$, see Section~8 in \cite{Ko2}; the rest of the Jordan 
structure remains the same; 

3) $c_{j_0}''$ belongs to the closure of $c_{j_0}$ (eventually, 
$c_{j_0}''=c_{j_0}$).  

When passing from $c_{j_0}'$ to $c_{j_0}''$ the index of rigidity decreases 
by at least $2$. If the change 2) can take place by changing the JNF of 
the restriction of $A_{j_0}$ or $M_{j_0}$ to some diagonal block, then we 
perform this change and further the proof is done as in the case of 
diagonalizable matrices.

$9^0$. If for the change 2) one has to change a block above the diagonal, 
and if there are at least $3$ diagonal blocks, then 
one proceeds as in $5^0$ and 
one proves that $d''<d'$ exactly in the same way. 

Indeed, at the first step 
one replaces two diagonal blocks (defining irreducible representations) 
by a single one (defining their semi-direct sum). Namely, using 
Notation~\ref{semidirect}, one chooses the block $Q_{j_0}$ such that 
the change 2) to take place. Then one chooses the block $Q_1$ such that 
condition (\ref{A_j}) or (\ref{M_j}) to hold (recall that $A_1$ and $M_1$ 
are with distinct eigenvalues, therefore changing the block $Q_1$ 
while keeping $P_1$ and $R_1$ the same does not change the conjugacy class of 
$A_1$ or $M_1$). 

The next steps are as in $5^0$.   

$10^0$. If there are just two diagonal blocks, not from Case B),  
then one first 
constructs a block upper-triangular tuple (with trivial centralizer) 
defining a semi-direct sum of the 
representations defined by the diagonal blocks but without changing the 
class $c_{j_0}'$. 

Then conjugate the tuple with a block upper-triangular 
matrix so that the matrix $A_{j_0}$ or $M_{j_0}$ 
to be in JNF (hence, it will be block diagonal as well). 
After this perform a change 
$A_{j_0}\mapsto A_{j_0}+\varepsilon U$ or 
$M_{j_0}\mapsto M_{j_0}+\varepsilon U$, $\varepsilon \in ({\bf C},0)$ 
where only the left lower block of $U$ is 
non-zero and is not of the form $R_{j_0}X-XP_{j_0}$; 
$U$ is chosen such that for $\varepsilon \neq 0$ one has 
$A_{j_0}\in c_{j_0}''$ 
(resp. $M_{j_0}\in C_{j_0}''$). 

To preserve condition (\ref{A_j}) or (\ref{M_j}) one looks then for 
deformations of the matrices $A_j$ or $M_j$, $j\neq j_0$, analytic in 
$\varepsilon$. Such a deformation exists, see the description of the 
``basic technical tool'' in \cite{Ko2} (one conjugates the matrices 
$A_j$ or $M_j$, $j\neq 0$, with matrices which are analytic deformations of 
$I$). 

\begin{lm}\label{Burnside}
For $\varepsilon \neq 0$ small enough the constructed tuple is irreducible.
\end{lm}

The lemma implies the theorem in this case.

$11^0$. If the two diagonal blocks are from Case B), then one change 2) is not 
sufficient to make the index of rigidity $\leq 0$. Hence, at least two 
changes are necessary. With the first of them we construct the semi-direct 
sum of representations defined by the two diagonal blocks; this time we change 
one of the JNFs for $j=j_*>1$. When performing this change we change the 
block $Q_{j_*}$ and then we change $Q_1$ to restore condition 
(\ref{A_j}) or (\ref{M_j}). 

Suppose that the second change must take place for $j=j_0\neq j_*$. 
Then after the second change 2) 
(performed as in $10^0$, using an analytic deformation) one has an 
irreducible representation by full analogy with Lemma~\ref{Burnside}.

If $j_*=j_0$ (and, say, $j_0=2$), then there are two 
possibilities. Either this JNF has a single eigenvalue, or it 
is with two double eigenvalues and three Jordan blocks. 
In the first case one can assume that the couple 
$A_2$, $U$ (resp. $M_2$, $U$) looks like this (after the analog of the 
conjugation from $10^0$): 

\[ A_2=\left( \begin{array}{cccc}a&1&0&0\\0&a&0&\underline{1}\\0&0&a&1\\0&0&0&a
\end{array}\right) ~~~,~~~U=\left( \begin{array}{cccc}0&0&0&0\\
0&0&0&0\\1&0&0&0\\0&0&0&0\end{array}\right) ~.\]
We underline the unit which is introduced after the first change 2). Its 
introduction results in changing the JNF like this: 
$\{ 2,2\} \mapsto \{ 3,1\}$. 
In the second case the couple looks like this: 

\[ A_2=\left( \begin{array}{cccc}a&0&\underline{1}&0\\0&b&0&0\\0&0&a&0\\0&0&0&b
\end{array}\right) ~~~,~~~U=\left( \begin{array}{cccc}0&0&0&0\\
0&0&0&0\\0&0&0&0\\0&1&0&0\end{array}\right) ~.\]

For the rest the proof is carried out as in $8^0$ -- $10^0$. 
The theorem is proved.~~~~~$\Box$

\subsection{Proofs of the lemmas\protect\label{prlm}}

{\bf Proof of Lemma~\ref{lmdim}:}\\ 

To obtain $d''$ one must add $l(n-l)$ to $d'''$, the dimension of the 
variety of block 
upper-triangular tuples as in the lemma (truly block upper-triangular, 
not only up to conjugacy). Indeed, $l(n-l)$ is the size of the left lower 
block and adding this corresponds to taking into account the possibility 
to conjugate such a tuple by matrices of the form 
$\left( \begin{array}{cc}I&0\\X&I\end{array}\right)$. 

One has $d'''=\Delta _1+\Delta _2+\Delta _3$ where 
$\Delta _1=\sum _{j=1}^{p+1}d_j^1-l^2+1$, 
$\Delta _2=\sum _{j=1}^{p+1}d_j^2-(n-l)^2+1$ and 
$\Delta _3=\sum _{j=1}^{p+1}s_j-l(n-l)$ (the contributions to $d'''$ from the 
$P$-, $R$- and $Q$-block). 

On the other hand, 
$d_j=d_j^1+d_j^2+2s_j$ (this can be deduced from Remark~\ref{d_jdiag}). Hence, 
$d'''=\sum _{j=1}^{p+1}d_j-\sum _{j=1}^{p+1}s_j-n^2+l(n-l)+2=
\sum _{j=1}^{p+1}d_j-\delta -n^2-l(n-l)+2$ and 
$d''=\sum _{j=1}^{p+1}d_j-\delta -n^2+2$. 

One has $d'=\sum _{j=1}^{p+1}d_j-n^2+1=d''+\delta -1$. Hence, for $\delta >1$ 
one has $d'>d''$.~~~~~$\Box$\\

{\bf Proof of Lemma~\ref{deform}:}\\ 

We transform the proof of the lemma into finding the cases when 
$\delta \leq 1$. 
%By Lemma~\ref{lmdim}, in all other cases the DSP is 
%solvable for the classes $C(P_j)\times C(R_j)$. 

\begin{stat}\label{stat1}
One has~~~$s_j\geq r_j^1(n-l)~~~(A)$~~~and~~~$s_j\geq r_j^2l~~~(B)$~~~(see 
Notation~\ref{semidirect}).
\end{stat}

{\em Proof:}\\  

Use Remark~\ref{semidirectrem} (and the 
notation from it) and Lemma~\ref{Ext}. Denote by 
$\mu '$ (resp. $\mu ''$) the biggest among the numbers 
$m_j'$ (resp. $m_j''$).  
Then $s_j\geq l(n-l)-\mu '(n-l)=r_j^1(n-l)$ because 
$\sum _{i=1}^sm_i'm_i''\leq \mu '\sum _{i=1}^sm_i''=\mu '(n-l)$. In the same 
way $s_j\geq l(n-l)-\mu ''l=r_j^2l$.~~~~~$\Box$  

\begin{rem}\label{equ}
Inequality $(A)$ becomes an equality exactly if $m_i''=0$ whenever 
$m_i'<\mu '$. Inequality $(B)$ becomes an equality exactly if 
$m_i'=0$ whenever $m_i''<\mu ''$. 
\end{rem}

\begin{stat}\label{stat2}
If for some index $j>1$ (say, $j=2$) one has $r_j^1=0$, $r_j^2>0$, then 
one has $\delta \geq 2$. The same is true if $r_j^1=r_j^2=0$ and $c_j$ is 
not scalar. The same is true if $r_j^1>0$, $r_j^2=0$.
\end{stat}

{\em Proof:}\\ 

Consider the first and the second of the three claims. By $(A)$ one has 
$s_3+\ldots +s_{p+1}\geq (r_3^1+\ldots +r_{p+1}^1)(n-l)\geq l(n-l)$; 
recall that $s_1=l(n-l)$. In the first claim one has 
also $s_2\geq r_j^2l\geq 2$, hence, $\delta \geq 2$. In the second claim 
the conjugacy class $c_2$ defines the MV $(l,n-l)$ and one 
has $s_2=l(n-l)\geq 2$ and again $\delta \geq 2$. The third claim is 
proved in the same way as the first one using $(B)$.~~~~~$\Box$\\ 

\begin{conv}
From now till the end of the proof of the lemma 
we assume (using the above statement) 
that for all indices $j>1$ one has $r_j^1>0$, $r_j^2>0$.
\end{conv}

\begin{stat}\label{stat3}
If $p\geq 3$, then $\delta \geq 2$.
\end{stat} 

{\em Proof:}\\

It suffices to consider the following two cases 
(up to permutation of the indices $j>1$):

1) $r_2^1\geq l/2$, $r_3^1\geq l/2$, $r_4^1>0$;

2) $r_2^1>0$, $l/2>r_j^1>0$ for $j>2$. 

In case 1) one has $s_2+s_3\geq l(n-l)$ (see $(A)$), 
$s_4\geq n-l\geq 2$, $s_1=l(n-l)$, so $\delta \geq 2$, see Lemma~\ref{Ext}.

In case 2) recall first that $r_j^2>0$ for $j>2$. 
For $j=3,4,\ldots ,p+1$ one has $s_j>r_j^1(n-l)$, 
i.e. $s_j\geq r_j^1(n-l)+1$, see Statement~\ref{stat1} 
and Remark~\ref{equ}. One has $s_1=l(n-l)$, 
$s_2\geq r_2^1(n-l)$ (see $(A)$), hence,  
$s_1+\ldots +s_{p+1}\geq 2l(n-l)+2$ and again $\delta \geq 2$.~~~~~$\Box$

\begin{conv}
From now till the end of the proof of the lemma we assume that $p=2$, 
see Statement~\ref{stat3}.
\end{conv}

\begin{stat}\label{stat4}
If $r_2^1+r_3^1\geq l+1$ or $r_2^2+r_3^2\geq n-l+1$, then $\delta \geq 2$.
\end{stat}

Indeed,  if $r_2^1+r_3^1\geq l+1$, then (see $(A)$) 
$s_2+s_3\geq (l+1)(n-l)\geq l(n-l)+2$ 
and $\delta \geq 2$. In the same way if $r_2^2+r_3^2\geq n-l+1$, then 
$s_2+s_3\geq l(n-l+1)\geq l(n-l)+2$ and $\delta \geq 2$.~~~~~$\Box$

\begin{stat}\label{stat5}
If $l$ is even and $r_2^1=r_3^1=l/2$, then $\delta \geq 2$, 
except in Cases B), C) and F) from the lemma.
\end{stat}

{\em Proof:}\\ 

$1^0$. If $l=2$, then $n-l=2$ and one has $\delta \leq 1$ only 
in one of Cases B) or C) from the lemma.\\ 

$2^0$. If $l\geq 4$, then $\delta \geq 2$. 
Indeed, to avoid Case A) 
from Lemma~\ref{betaequalalpha} for the block $P$, 
one must suppose that at least one of 
the two matrices $P_2$ and $P_3$ (say, $P_2$) has at least 
three distinct eigenvalues. Assume that the MV of $P_2$ looks like this: 
$(m_1',\ldots ,m_s')$, with $m_1'=\mu '>m_2'\geq \ldots \geq m_s'$ (the 
inequality $m_1'>m_2'$ results from $r_2^1=l/2$). 

If for at least two indices $i>1$ one has $m_i''\neq 0$, then 
for them one has $m_i'm_i''<\mu 'm_i''$ and 
$\sum _{i=1}^sm_i'm_i''\leq \mu '\sum _{i=1}^sm_i''-2=\mu '(n-l)-2$. Hence, 
$s_2\geq r_2^1(n-l)+2$ (see Remark~\ref{semidirectrem}), 
$s_3\geq r_3^1(n-l)$ (see $(A)$) and $\delta \geq 2$.\\ 

$3^0$. If for only one index $i>1$ one has $m_i''\neq 0$ (i.e. $R_2$ 
has only two different eigenvalues), then similarly 
$s_2\geq r_2^1(n-l)+1$ with equality only if the MV of $P_2$ equals 
$(l/2,l/2-1,1)$ and the two eigenvalues, of the two greatest multiplicities, 
are eigenvalues of $R_2$ as well; moreover, its only eigenvalues.\\  

$4^0$. If $P_3$ has at least three different eigenvalues, 
then in the same way $s_3\geq r_3^1(n-l)+1$ and, hence, $\delta \geq 2$. 
So the only possibility to have $\delta \leq 1$ is the MV of $P_3$ to be 
$(l/2,l/2)$. 
If $n-l=2$, then $\delta \leq 1$ only in Case F). If $n-l>2$, then $R_3$ must 
have at least three distinct eigenvalues (otherwise condition 
$(\alpha _{n-l})$ fails for the block $R$) and $s_3\geq lr_3^2+1$. One has 
also $s_2\geq lr_2^2+1$ (to be checked directly), 
hence, again $\delta \geq 2$.~~~~~$\Box$

\begin{stat}\label{stat6}
Suppose that $r_2^1+r_3^1=l$. 
If $r_2^1>l/2$, $r_3^1<l/2$ or $r_2^1<l/2$, $r_3^1>l/2$, then $\delta \geq 2$ 
except in Cases D), E) from the lemma.
\end{stat}

{\em Proof:}\\ 

$1^0$. Without loss of generality we assume that $r_2^1>l/2$, $r_3^1<l/2$. 
If $l=3$ and $n-l=2$ or $n-l=3$, 
then one has $\delta \leq 1$ only 
in Case E) with $q=1$ or in Case D) of the lemma. Indeed, $s_j$ is minimal 
only if all eigenvalues of $R_j$ are 
eigenvalues of $P_j$ as well for $j=2,3$.\\  

$2^0$. If $l\geq 5$ and $n-l\geq 4$, then $\delta \geq 2$. 
Indeed, if the MVs 
of $P_3$ and $R_3$ equal respectively $(m_1',\ldots ,m_s')$, 
$(m_1'',\ldots ,m_s'')$, with 
$m_1'=\mu '>m_2'\geq \ldots \geq m_s'$, then one has 
$m_i'm_i''\leq (\mu '-1)m_i''$ for $i>1$ and $m_i''>0$; hence, 

\[ \sum _{i=1}^sm_i'm_i''\leq \mu 'm_1''+(\mu '-1)\sum _{i=2}^sm_i''
=\mu '\sum _{i=1}^sm_i''-\sum _{i=2}^sm_i''=\mu '(n-l)-\sum _{i=2}^sm_i''~.\] 
If $\sum _{i=2}^sm_i''\geq 2$, then $s_3\geq r_3^1(n-l)+2$ 
(see Remark~\ref{semidirectrem}), 
$s_2\geq r_2^1(n-l)$ (see $(A)$) and 
$\delta \geq 2$. So $\delta$ can be $\leq 1$ only in case that 
$\sum _{i=2}^sm_i''=1$, i.e. the MV of $R_3$ is of the form $(n-l-1,1)$. 
If this is so, then the MV of $R_2$ is $(1,\ldots ,1)$ (otherwise 
$(\alpha _{n-l})$ fails for the block $R$), i.e. $R_2$ has 
distinct eigenvalues. Hence, $s_2\geq l(n-l-1)$ whatever the eigenvalues of 
$P_2$ are.  

But then $s_3$ is minimal if and only if the MV of $P_3$ equals 
$(l-1,1)$ and $P_3$ has the same eigenvalues as $R_3$ (the proof of this 
is left for the reader). In this case $s_3=(l-1)+(n-l-1)=n-2$, hence, 
$s_2+s_3\geq l(n-l)+n-l-2$ and for $n-l\geq 4$ one has $\delta \geq 2$.\\ 

$3^0$. If $l\geq 5$ and 
$n-l=2$, and if $P_2$ has at least $4$ distinct 
eigenvalues, then $s_2\geq l+2$. Indeed, $s_2$ is minimal only if each 
eigenvalue of $R_2$ is eigenvalue of $P_2$ as well. 

In such a case one has 
$s_2=2l-m_{i_1}'-m_{i_2}'$ where $m_{i_1}'$, $m_{i_2}'$ are the 
multiplicities of the eigenvalues of $R_2$ as eigenvalues of $P_2$. As 
$m_{i_1}'+m_{i_2}'\leq l-2$ (there are at least two more eigenvalues of $P_2$, 
each of multiplicity $\geq 1$), one gets $s_2\geq l+2$. In a similar way, 
$s_3\geq l$, with equality when $P_3$ has two eigenvalues which are 
eigenvalues of $R_3$ as well, hence, $\delta \geq 2$. 

If $P_2$ has exactly three distinct eigenvalues, then one has $s_2\geq l+1$ 
with equality exactly if the eigenvalue which is not eigenvalue of $P_2$ is 
simple. Hence, $\delta \leq 1$ only in Case~E) from the lemma.\\  

$4^0$. If $l\geq 5$ and $n-l=3$, then at least one of the matrices 
$R_2$, $R_3$ must have $3$ distinct eigenvalues (otherwise $(\beta _3)$ fails 
for the block $R$). The respective quantity 
$s_j$ must be $\geq 2l=r_j^2l$, see $(B)$. 
If the other matrix $R_j$ ($j=2$ or $3$) 
has also $3$ distinct eigenvalues, 
then $s_2+s_3\geq 4l>3l+2$ and $\delta \geq 2$.  

If the MV of the other matrix $R_j$ (say, $R_3$) equals $(2,1)$, then $s_3$ 
is minimal exactly if $P_3$ has the same eigenvalues as $R_3$, of 
multiplicities $l-1$ and $1$. In this case $s_3=l+1$. But then $P_2$ must be 
with distinct eigenvalues (otherwise $(\alpha _l)$ fails for 
the block $P$), $s_2\geq 3l-3$, and $\delta >2$.\\ 

$5^0$. If $l=4$, then one can have $r_2^1>2$, $r_3^1<2$ only if $P_3$ has four 
distinct eigenvalues and the MV of $P_3$ is $(1,3)$. We let the reader 
check oneself that in all possible cases ($n-l=2,3$ or $4$) one has 
$\delta \geq 2$.~~~~~$\Box$\\  
    
The lemma follows from Statements~\ref{stat2}, \ref{stat3}, \ref{stat4}, 
\ref{stat5} and \ref{stat6}.~~~~~$\Box$\\ 

{\bf Proof of Lemma~\ref{Burnside}:}\\ 

Denote by ${\cal T}$, 
the matrix algebra of all block upper-triangular matrices 
with square diagonal blocks of sizes $l$ and $n-l$. A priori the 
representation defined by the deformed matrices is either 
irreducible (and the corresponding matrix algebra is $gl(n,{\bf C})$) or 
is reducible and defines a matrix algebra which up to analytic 
conjugation equals ${\cal T}$ (the statement results from a more 
general one which can be found in \cite{Ko7}). 
The second case, however, is impossible 
because such a conjugation of $A_{j_0}$ or $M_{j_0}$ (with a matrix 
$I+O(\varepsilon )$) cannot make the 
left lower block of $U$ disappear (because it is not of the form 
$R_{j_0}X-XP_{j_0}$).~~~~~$\Box$

\section{Proof of Theorem~\protect\ref{n+1}\protect\label{proofofn+1}}
\subsection{Proof in the case of matrices $A_j$\protect\label{caseofA_j}}

\begin{defi}\label{regulardefi}
A conjugacy class is called {\em regular} if to every eigenvalue there 
corresponds a single Jordan block of size equal to the multiplicity of the 
eigenvalue.
\end{defi}

\begin{rem}\label{regularrem}
The JNFs of all regular conjugacy classes correspond to each other (see 
Definition~\ref{corrdefi}) and, in particular, to the diagonal JNF with 
distinct eigenvalues and to the JNF with a single eigenvalue and a single 
Jordan block of size $n$.
\end{rem}

\begin{prop}\label{reg2n}
The DSP is positively solvable for classes $c_j$ where $c_1$ is regular and 
one has $r_2+\ldots +r_{p+1}\geq n+1$.
\end{prop}

The proposition implies the theorem in the case of matrices $A_j$. 
To prove the proposition we need the following lemma.

\begin{lm}\label{nilp2n}
The DSP is positively solved for tuples of nilpotent conjugacy classes $c_j$ 
with $r_1+\ldots +r_{p+1}\geq 2n$ in which $r_1=n-1$, i.e. the conjugacy class 
$c_1$ has a single Jordan block of size $n$.
\end{lm}

The lemma is a particular case of the results in \cite{Ko8}. It follows also 
from the ones in \cite{C-B}.\\ 

{\em Proof of the proposition:}\\ 

Given an irreducible tuple of nilpotent matrices $A_j$ satisfying the 
conditions of the lemma one can deform it analytically into an irreducible 
tuple of 
matrices $A_j'$ where for each $j$ either $J(A_j')=J(A_j)$ or $J(A_j')$ 
corresponds to $J(A_j)$. The eigenvalues of the matrices $A_j'$ must be close 
to $0$. These statements can be deduced from \cite{Ko2}, see the definition 
of the basic technical tool there which is a way to deform analytically 
tuples of matrices with trivial centralizers; compare also with Lemma~53 
from \cite{Ko2}. 

Thus one obtains the positive solvability of the DSP for all tuples of 
JNFs $J(c_j)$ satisfying the condition $r_2+\ldots +r_{p+1}\geq n+1$; 
see Definition~\ref{corrdefi} and Remarks~\ref{corrrems} (especially part 
2) of them). However, solvability is proved only for eigenvalues close to $0$. 

By multiplying the tuples of matrices 
$A_j'$ by non-zero complex numbers (i.e. $(A_1',\ldots ,A_{p+1}')\mapsto 
(gA_1',\ldots ,gA_{p+1}'), g\in {\bf C}^*$) one can obtain 
irreducible tuples with the same JNFs as $A_j'$ 
and with any eigenvalues whose sum 
(taking into account the multiplicities) is $0$. 
This proves the proposition.~~~~~$\Box$

\subsection{Proof for matrices $M_j$}

Suppose that for some conjugacy classes $C_j$ satisfying the conditions of the 
theorem there exist no irreducible tuples. Then there exist tuples with 
trivial centralizers. This follows from Theorem~\ref{weakDSPdiag} and from 
Lemma~\ref{betaequalalpha}. 

Each such tuple can be conjugated to a block 
upper-triangular form in which the diagonal blocks define irreducible or 
one-dimensional representations. Denote by $s_1$, $\ldots$, $s_{\nu}$ the 
sizes of the diagonal blocks. We say that these sizes (considered up to 
permutation) define the {\em type} of the tuple. 
The tuple is called {\em maximal} 
if there is no tuple with trivial 
centralizer and of type 
$s_1'$, $\ldots$, $s_h'$ such that $h<\nu$ and the sizes $s_i'$ are 
obtained from the sizes $s_j$ by one or several operations of the form 
$(s_{j_1},s_{j_2})\mapsto s_{j_1}+s_{j_2}$. We say that the type 
$s_1'$, $\ldots$, $s_h'$ is {\em greater} than the type 
$s_1$, $\ldots$, $s_{\nu}$. 

\begin{lm}\label{constructtuple}
Given a maximal tuple of matrices $M_j$ one can construct a tuple 
of matrices $A_j\in c_j$ of the same 
type, with trivial 
centralizer, with $M_j=\exp (2\pi iA_j)$ (up to conjugacy) where for $j>1$ 
the matrix $A_j$ has no couple of eigenvalues whose difference is a 
non-zero integer.
\end{lm}

The lemma is proved in the next subsection.

\begin{rem}
The condition ``$M_j=\exp (2\pi iA_j)$ (up to conjugacy)'' is introduced with 
the aim to use  
the fact that the monodromy operators of the Fuchsian system 
$dX/dt=(\sum _{j=1}^{p+1}A_j/(t-a_j))X~(**)$ 
in the absence of non-zero integer 
differences between the eigenvalues of the matrices $A_j$ 
equal (up to conjugacy)
$\exp (2\pi iA_j)$. See the definition of the monodromy 
operators in the Introduction of \cite{Ko2}.  
\end{rem}

For the tuple of matrices $A_j$ from the lemma one has that they can be 
analytically deformed into an irreducible tuple of such matrices. 
Indeed, for their conjugacy classes the DSP is positively solved 
(this is already proved in Subsection~\ref{caseofA_j}) and all 
reducible tuples 
from these classes belong to the closure of the variety of irreducible tuples, 
see Theorem~6 from \cite{Ko3}. 

All irreducible tuples of matrices $A_j^0$ close to tuples $A_j$ 
from the lemma 
define Fuchsian systems $dX/dt=(\sum _{j=1}^{p+1}A_j^0/(t-a_j))X~(***)$ 
whose monodromy groups 
must be (up to conjugacy) 
from the type of the tuple of matrices $M_j$ from the lemma. This follows from 
the tuple of matrices $M_j$ being maximal. 

Consider the monodromy operators (denoted also by $M_j$) of systems $(**)$ 
with matrices-residua $A_j$. One has $M_j=\exp (2\pi iA_j)$ (up to conjugacy) 
and there is a bijection between the eigenvalues of the matrices $A_j$ 
and the ones of the matrices $M_j$. For each diagonal block the sum of the 
eigenvalues of the matrices $A_j$ from the lemma is $0$. Hence, the sum of the 
same eigenvalues of the matrices $A_j^0$ is also $0$. 
If the monodromy group of system $(***)$ is of the type of the one of system 
$(**)$, then by Theorem~5.1.2 from \cite{Bo} it should be possible to 
conjugate the tuple of matrices $A_j^0$ to a block upper-triangular form 
with blocks as in the type of the matrices $M_j$. This contradicts the 
irreducibility of the tuple of matrices $A_j^0$.

\begin{rem}
When applying Theorem~5.1.2 from \cite{Bo} we use the fact that there are 
no non-zero integer differences between the eigenvalues of the matrices 
$A_j$. Thus to each eigenvalue $\sigma$ of $M_j$ of a given multiplicity there 
corresponds only one eigenvalue $\lambda$ of $A_j$ (which is 
of the same multiplicity) 
where $\sigma =\exp (2\pi i\lambda )$. 
Theorem~5.1.2 from \cite{Bo} speaks about the exponents (i.e. the eigenvalues 
of the matrices $A_j$) corresponding to an invariant subspace. In the 
absence of non-zero integer differences these exponents are defined by 
the eigenvalues of the monodromy operators in a unique way.
\end{rem}

The theorem is proved.~~~~~$\Box$

\subsection{Proof of Lemma~\protect\ref{constructtuple}}

$1^0$. One can construct for each size $s_i$ of the type a tuple of 
matrices $A_{i,j}^*$ such that one has 
(up to conjugacy) $\exp (2\pi i A_{i,j}^*)=M_{i,j}^*~(*)$   
where $M_{i,j}^*$ are the 
restrictions of the matrices $M_j$ to the diagonal block of size $s_i$, and 
the matrices $A_{i,j}^*$ define an irreducible or 
one-dimensional representation. In the one-dimensional case the claim is 
evident. In the irreducible case one can construct a Fuchsian system with 
matrices-residua equal up to conjugacy to $A_{i,j}^*$ 
(where $A_{i,j}^*$ satisfy $(*)$) the real parts of whose eigenvalues 
can be chosen to 
belong to $[0,1)$ for $j>1$ (to avoid non-zero integer differences between 
eigenvalues); the construction is explained in \cite{ArIl}.\\  

$2^0$. Consider the tuple of matrices $A_j'$ which are block-diagonal their 
restrictions to each diagonal block of size $s_i$ being equal to 
the blocks $A_{i,j}^*$ from $1^0$. We complete them (in $3^0$) 
by adding 
entries in the blocks above the diagonal (the newly obtained 
matrices are denoted by $A_j$) so that one would have 
$\exp (2\pi i A_j)=M_j$ up to conjugacy. We do this for $j>1$ and then we 
define $A_1$ so that $A_1+\ldots +A_{p+1}=0$. As $A_1$ has distinct 
eigenvalues, whatever entries we add in the blocks above the diagonal, 
they do not change the conjugacy class of $A_1$. 
As $\exp (2\pi i A_1')=M_1$ up to conjugacy, one will also 
have $\exp (2\pi i A_1)=M_1$ up to conjugacy.\\  

$3^0$. One can conjugate the matrix 
$M_j$ by a block upper-triangular matrix $B_j$ 
so that 
the diagonal blocks of $(B_j)^{-1}M_jB_j$ of sizes $s_i$ 
to be in JNF and in the blocks above the 
diagonal non-zero entries to be present only in positions $(i,j)$ such that 
the $i$-th and $j$-the eigenvalues coincide. For each eigenvalue 
$\sigma _{k,j}$ of $M_j$ denote by $M_j(\sigma _{k,j})$ the matrix 
whose restriction to the rows and columns of the eigenvalue $\sigma _{k,j}$ 
are the same as the ones of $(B_j)^{-1}M_jB_j$ and the rest of its entries 
are $0$.

One can conjugate the matrices $A_j'$ by block-diagonal matrices $D_j$ so that 
the matrix $(D_j)^{-1}A_j'D_j$ to be in JNF and for each diagonal block 
there to hold $\exp (2\pi iA_{i,j}^*)=M_{i,j}^*$ (up to conjugacy). 

Set $(D_j)^{-1}A_j'D_j=\sum _{k,j}\lambda _{k,j}A_j'(\lambda _{k,j})$ where 
$\lambda _{k,j}$ are the distinct eigenvalues of $A_j'$ and 
$A_j'(\lambda _{k,j})$ is the matrix 
whose restriction to the rows and columns of the eigenvalue $\lambda _{k,j}$ 
is the same as the one of $(D_j)^{-1}A_j'D_j$ and the rest of its entries 
are $0$. Define the matrices $A_j(\lambda _{k,j})$ by analogy with the 
matrices $A_j'(\lambda _{k,j})$.

Recall that one has $\sigma _{k,j}=\exp (2\pi i \lambda _{k,j})$. 
Hence, for each diagonal block and for each couple $(k,j)$ the restrictions of 
the matrices $A_j'(\lambda _{k,j})-\lambda _{k,j}I$ and 
$M_j(\sigma _{k,j})-\sigma _{k,j}I$ to it are equal. 

Define the matrices 
$(D_j)^{-1}A_jD_j$ by the rule for all $(k,j)$ the 
matrices $A_j(\lambda _{k,j})-\lambda _{k,j}I$ and 
$M_j(\sigma _{k,j})-\sigma _{k,j}I$  
to be equal. The rule implies that the JNFs of the matrices 
$(B_j)^{-1}M_jB_j$ and $(D_j)^{-1}A_jD_j$, hence, of $M_j$ and $A_j$, 
coincide. As there are no non-zero integer differences between eigenvalues of 
$A_j$, one has also $\exp (2\pi iA_j)=M_j$ (up to conjugacy).\\ 
 
$4^0$. The tuple of matrices $A_j$ thus constructed might fail to be 
with trivial centralizer. Hence, the tuple must define a direct sum of 
representations (this follows from $A_1$ being with distinct eigenvalues). 
So conjugate it to a block-diagonal form where each block 
(we call these blocks {\em big blocks}) is 
small-block upper-triangular and with trivial centralizer. The small blocks 
are of sizes $s_i$. 

As in Lemma 24 from \cite{Ko3} one shows 
that if there are two big blocks of sizes $u,v$ where $u\geq 3$, $v\geq 2$, 
then one can deform the tuple into one in which these two big blocks are 
replaced by a single big block of size $u+v$ (with trivial centralizer 
and with the same small blocks as the two big blocks) 
while the other big blocks remain the same. The statement holds also if 
$u=v=2$, $p=2$ (see again Lemma 24 from \cite{Ko3}) 
and for at least one index $j\geq 2$ the 
restrictions of the tuple to 
the two big blocks belong to different conjugacy classes, or if $u=v=2$, 
$p\geq 3$ and no matrix is scalar.   

If there is a big block $B$ of size $1$, then it follows from 
$r_2+\ldots +r_{p+1}\geq n+1$ that for at least one of the other big blocks 
$B'$ one has Ext$^1(B,B')\geq 1$. Indeed, without loss of 
generality one can assume that the restrictions of the matrices to the 
block $B$ equal $0$ for all values of $j$. Hence, for each other big block 
$B'$ one has Ext$^1(B,B')=\rho (B')-2\sigma (B')$ where 
$\rho (B')$ is the 
sum of the ranks $r_j(B')$ of the matrices $A_j|_{B'}$ and $\sigma (B')$ 
is the size of 
$B'$. (One subtracts $\sigma (B')$ once because the sum of the 
matrices $A_j$ is $0$ and once to factor out conjugation 
with block upper-triangular matrices; see the proof of Lemma~\ref{Ext}.) 

If for all blocks $B'$ one has Ext$^1(B,B')\leq 0$, then 
one has   

\[ 0\geq \sum _{B'}(\rho (B')-2\sigma (B'))=
(\sum _{j=1}^{p+1}\sum _{B'}r_j(B'))-2(n-1)\geq (\sum _{j=1}^{p+1}r_j)-2(n-1)\]
i.e. $\sum _{j=2}^{p+1}r_j\leq n-1$ (recall that $r_1=n-1$) 
which is a contradiction.
 
Hence, one can replace the two 
blocks $B$, $B'$ by a single big block of size $\sigma (B')+1$. 

There remains to be considered the case when there is no big block of size 
$1$ or $\geq 3$, i.e. all big blocks are of size $2$; moreover, $p=2$, and 
for $j>1$ the 
restrictions of the matrices $A_j$ to the big blocks belong to one and the 
same conjugacy class. In this case one has $r_2+r_3=n$, i.e. the case has 
not to be considered.~~~~~$\Box$

Author's address: Universit\'e de Nice -- Sophia Antipolis, Laboratoire de 
Math\'ematiques, Parc Valrose, 
06108 Nice, Cedex 2, France; e-mail: kostov@math.unice.fr

\end{document}